\documentclass[12pt]{amsart}
\usepackage[a4paper,hmargin=2cm,vmargin=3cm]{geometry}
\usepackage[T1]{fontenc}
\usepackage[utf8]{inputenc}
\usepackage{color}
\usepackage{amsmath}
\usepackage{amsthm}
\usepackage{amssymb}
\usepackage{graphicx}
\usepackage[style=alphabetic,maxnames=10,giveninits=true,backend=biber,doi=false,url=false]{biblatex}
\usepackage[shortlabels]{enumitem}
\usepackage{hyperref}
\usepackage{amsaddr}
\usepackage{ulem}
\definecolor{notefontcolor}{rgb}{0.800781, 0.800781, 0.800781}
\definecolor{grey30}{rgb}{0.7,0.7,0.7}

\hypersetup{unicode=true,
  colorlinks=true, citecolor=blue, linkcolor=blue,
  pdftitle=Complexity of local maxima of given radial derivative for mixed p-spin Hamiltonians,
  pdfauthor={D. Belius, M. Schmidt}}

\DeclareSourcemap{
  \maps[datatype=bibtex]{
    \map{
      \step[fieldset=issn, null]
      }
    }
  }
\numberwithin{equation}{section}

\theoremstyle{plain}
\newtheorem{theorem}{Theorem}[section]
\newtheorem{lemma}[theorem]{Lemma}
\newtheorem{conjecture}[theorem]{Conjecture}

\newtheorem{rem}[theorem]{Remark}

\theoremstyle{remark}

\DeclareMathOperator{\Cov}{Cov}
\newcommand{\beq}{ \begin{equation}} 
\newcommand{\eeq}{\end{equation}} 
\newcommand{\bea}{\begin{aligned}}
\newcommand{\eea}{\end{aligned}}

\newcommand{\Num}{{\mathcal{N}}}
\newcommand{\R}{{\mathbb{R}}}
\renewcommand{\P}{{\mathbb P}}
\newcommand{\E}{\mathbb E}

\newcommand{\sph}{{\mathrm{sp}}}

\newcommand{\xip}{\boldsymbol{\xi}'}
\newcommand{\xipp}{\boldsymbol{\xi}''}
\newcommand{\xipa}{{\boldsymbol{\xi}'_{\alpha}}}
\newcommand{\xippa}{{\boldsymbol{\xi}''_{\alpha}}}
\newcommand{\const}{\mathrm{const}}
\definecolor{mypink1}{rgb}{0.858, 0.188, 0.478}

\begin{document}

\title[Complexity for given radial derivative]{Complexity of local maxima of given radial derivative for mixed $p$-spin Hamiltonians}

\author{David Belius, Marius A. Schmidt}
\address{Department of Mathematics and Computer Science, University of Basel, Switzerland}

\email{david.belius@cantab.net}
\email{m.schmidt@mathematik.uni-frankfurt.de}
\thanks{Both authors supported by SNSF grant 176918. Marius A. Schmidt supported by a DFG research grant, contract number 2337/1-1.}

\maketitle

\begin{abstract}
We study the number of local maxima with given radial derivative of spherical mixed $p$-spin models and prove that the second moment matches the square of the first moment on exponential scale for arbitrary mixtures and any radial derivative. This is surprising, since for the number of local maxima with given radial derivative and given energy the corresponding result is only true for specific mixtures \cite{subag2017complexity,arous2020geometry}.

We use standard Kac-Rice computations to derive formulas for the first and second moment at exponential scale, and then find a remarkable analytic argument that shows that the second moment formula is bounded by twice the first moment formula in this general setting.

This also leads to a new proof of a central inequality used to prove concentration of the number critical points of pure $p$-spin models of given energy in \cite{subag2017complexity} and removes the need for the computer assisted argument used in that paper for $3 \le p \le 10$.
\end{abstract}

\section{Introduction}
The spherical mixed $p$-spin models are a natural and general class of isotropic differentiable Gaussian random fields on the sphere. They are paradigmatic models of high-dimensional complex random landscapes and originate in spin glass theory \cite{sherrington1975solvable,derridaRandomEnergyModelLimit1980,grossSimplestSpinGlass1984,talagrandMultipleLevelsSymmetry2000,kosterlitzSphericalModelSpinGlass1976,crisantiSphericalpspinInteractionSpin1992,talagrandFreeEnergySpherical2006,crisantiSphericalSpinglassModel2004,mezard1987spin,talagrand2010mean,panchenko2013sherrington}. We study the number of local maxima of the field with fixed radial derivative. Specifically, we compute the second moment, and show that on an exponential scale it matches the first moment {\emph{for any mixed $p$-spin model and any radial derivative}}. This strongly suggests, but does not yet prove, that the number of local maxima of given radial derivative concentrates around its mean. Expressed in the spin glass terminology it thus strongly suggests that quenched and annealed complexity of local maxima of given radial derivative  {\emph{always}} coincides.

This is surprising, since for the previously studied number of local maxima (or critical points) at fixed radial derivative and fixed energy this is only true for all radial derivatives and energies for very special mixed $p$-spin models, namely the pure $p$-spin models and their perturbations \cite{auffinger2013random,AB13,subag2017complexity,arous2020geometry}.

To state our result, consider a mixed $p$-spin Hamiltonian $H_N$ with mixture $\xi$ on the unit sphere  $S_{N-1}\subset\R^N$, i.e $H_N$ is a centered Gaussian field on $S_{N-1}$ with covariance given by 
$$ \Cov[ H_N(\sigma), H_N(\tau)] = N \xi(\sigma \tau),$$ for $\sigma,\tau\in S_{N-1}$ and $\sigma\tau= \sum_{i=1}^{N}\sigma_i\tau_i$ the inner product. Any isotropic centered Gaussian field on the sphere must have a covariance of this form, and a $\xi$ gives a well-defined covariance for all $N$ if and only if it is of the form $\xi(x)= \sum_{p\geq 0} a_p x^p$ for $p\ge0$ with $\xi(1)<\infty$ \cite{schoenbergPositiveDefiniteFunctions1942a}. Our results apply to any $\xi$ of such a form with $a_0=a_1=1$ and radius of convergence greater than $1$, which is thus a very general class of isotropic Gaussian fields on the sphere. For these $\xi$ the field $H_N$ is almost surely smooth on a ball of radius larger than one (see Lemmata A.1 and A.7 \cite{TAPHTUB}).
For $\partial_r$ the derivative in radial direction the central object of our study is the number
$$ \Num(D) = \text{Number of local maxima of } H_N \mbox{ on } S_{N-1} \text{ with } \frac{1}{N}\partial_r H_N \in D.$$
Note that $H_N$ is defined on a ball so that we may speak of its radial derivative, but  we always consider local maxima with respect to the unit sphere. We will use what are by now standard Kac-Rice computations to show that (see Lemma \ref{lem: first mom})
\begin{equation}\label{eq: intro first mom} \frac{1}{N}\ln E[\Num(D)] \to \sup_{x\in D}I(x),
\end{equation}
where the function $I$ is given explicitly in terms of the shorthands $\xip=\xi'(1)$ and  $\xipp= \xi''(1)$ by
$$ I(x) := \frac{1}{2}\ln\left(\frac{\xipp}{\xip}\right)-\frac{x^2}{4\xipp}\frac{\xipp-\xip}{\xipp+\xip}+\Omega\left(\frac{x}{\sqrt{2\xipp}}\right),$$
for
\beq \label{eq: omegadef} \Omega(y):=\begin{cases}-\frac{1}{2}y\sqrt{y^2-2}+\ln\left(\frac{y+\sqrt{y^2-2}}{\sqrt{2}}\right)& \mbox{ for } y\geq \sqrt{2},\\ -\infty & \mbox{ else.}\end{cases} \eeq
Writing
\beq \label{def: rinfty} r_\infty := \inf\{ x\in \R : I(x)\geq 0\} = 2 \sqrt{\xipp}\quad \mbox{ and }\quad r_0 := \sup\{x\in\R: I(x)\geq 0\}\eeq
we note that $I$ is $-\infty$ below $r_\infty$ and strictly decreasing on $[r_\infty,\infty)$ with $I(r_0)=0$. This entails by Markov inequality that with high probability there are no local maxima of $H_N$ on $S_{N-1}$ with radial derivative significantly above $r_0$ or below $r_\infty$. 
Our main result is the following, which shows the second and first moment match on exponential scale for {\emph{any}} mixture $\xi$. 
\begin{theorem}[Matching moments on exponential scale]\label{thm: exp sec mom}
For all $x \in [r_\infty,r_0]$ we have 
\begin{equation}\label{eq: exp sec mom}
\lim_{\varepsilon \searrow 0}\lim_{N\rightarrow\infty}\frac{1}{N} \ln E[ \Num(\left[x-\varepsilon,x+\varepsilon\right])^2 ] - \frac{2}{N} \ln E[ \Num(\left[x-\varepsilon,x+\varepsilon\right]) ] = 0.
\end{equation}
Additionally  
\beq\label{eq: exp fst mom} \lim_{\varepsilon \searrow 0}\lim_{N\rightarrow\infty}\P(\Num(\left(-\infty,r_\infty-\varepsilon\right]\cup \left[r_0+\varepsilon,\infty\right))=0)=1.
\eeq
\end{theorem}
As mentioned above, it is surprising that the matching \eqref{eq: exp fst mom} of moments holds {\emph{for all}} $\xi$ with $a_0=a_1=0$ and $x\in[r_\infty,r_0]$. Proving the matching of the first and second moment on exponential scale is the main step when proving concentration around the mean using the second moment method. Since this main step is achieved by the our theorem we make the following conjecture.
\begin{conjecture}\label{conj}
For all $r\in [r_\infty,r_0]$ 
\begin{equation}\label{eq: quenched}
    \lim_{\varepsilon\searrow 0} \lim_{N\rightarrow\infty}\frac{1}{N}\ln \Num(\left[x-\varepsilon,x+\varepsilon\right]) =  I(x)
\end{equation}
in probability.
\end{conjecture} 
In particular, Conjecture \ref{conj} would imply that the total number of local maxima is in $\exp\left(N I(r_\infty)+o(N)\right)$ with high probability for all $\xi$.

One way to prove Conjecture \ref{conj} would be to show that the second moment is asymptotic to the first moment squared, i.e. a stronger version of \eqref{eq: exp sec mom}. This seems attainable using the methods of \cite{subag2017complexity} for mixtures with $a_2=0$, but is beyond the scope of this paper\footnote{The restriction $a_2=0$ is to have approximate independence of $M^{(1)}$ and $M^{(2)}$ for small $\alpha$, see proof of Lemma \ref{lem: sec mom}, \cite[Lemma 13 (4.8) and (2)]{subag2017complexity} and \cite[Lemma 14 (4.10) and (2)]{arous2020geometry}.}. Alternatively, if one could obtain a general concentration result for $\mathcal{N}$ such as is available e.g. for Lipschitz functions of independent Gaussian random variables then Conjecture \ref{conj} could be derived from the exponential scale matching of moments \eqref{eq: exp sec mom}.

To prove Theorem \ref{thm: exp sec mom} we follow \cite{subag2017complexity,arous2020geometry} by using the Kac-Rice formula to compute a variational formula for the second moment (see Lemma \ref{lem: sec mom}). The upshot is
$$ \frac{1}{N}\ln E[\Num(D)^2]\rightarrow \sup\left\{ S(x,\alpha):\alpha\in(-1,1), x\in D\right\},$$
where $S$ is given by \eqref{def: two-point-function} below.

It is easily checked that for all $x\in[r_\infty,r_0]$ one has $S(x,0)=2I(x)$. The main step of the proof of \eqref{eq: exp sec mom} is then an elementary but complicated computation that proves that $S(x,\alpha)$ is maximized at $\alpha=0$, which with \eqref{eq: intro first mom} gives \eqref{eq: exp sec mom}. Beyond the fact that this holds for \emph{any} mixture and $x\in[r_\infty,r_0]$, it is also pleasantly surprising that an analytic proof of it can be obtained.

For pure $p$-spin mixtures, i.e. for $\xi(x)=x^p$, a result equivalent to Theorem \ref{thm: exp sec mom} was proved in \cite{subag2017complexity}, which studied the number of critical points of pure $p$-spin models at fixed energy. The equivalence is due to the fact that at low energies most critical points are local maxima, and that for pure $p$-spin models $\partial_r H_N(\sigma) = p H_N(\sigma)$ for all $\sigma \in S_{N-1}$ almost surely and therefore a restriction on the radial derivative is equivalent to a restriction on the energy. Note however that for $3\le p \le 10$ the proof of \cite{subag2017complexity} is computer assisted, as a computer plot is employed to show that the formula corresponding to $S(x,\alpha)$ is maximized at $\alpha=0$ (see \cite[Proof of Lemma 7; Figure 1]{subag2017complexity}. Our proof on the other hand is fully analytic for all $\xi$.

Since its introduction as a tool for the mathematically rigorous study of spin glasses and mixed $p$-spin models \cite{fyodorov2004complexity,fyodorovCriticalBehaviorNumber2012,auffinger2013random}  the use of the Kac-Rice formulas has become standard. Without distinguishing between the different types of critical points that have been considered (all critical points, local maxima, saddles of fixed index, etc.), the first moment for pure mixtures is discussed in \cite{auffinger2013random,fyodorov2013high}, the second moment in \cite{subag2017complexity}, the first moment for mixed models in \cite{AB13,fyodorov2013high} and the second in \cite{arous2020geometry}. As mentioned above, for pure models \cite{subag2017complexity} shows that the complexity  concentrates around its mean. The results of \cite{arous2020geometry} extend this to mixed models that are small perturbation of pure models. The first moment in the presence of external field is computed in \cite{fyodorov2013high, beliusTrivialityGeometryMixed2022}.

Before going into details we discuss the structure of the paper. In Section \ref{sec: moments} we compute the limits in \eqref{eq: exp sec mom} in form of Lemma \ref{lem: first mom} for the first moment and Lemma \ref{lem: sec mom} for pairs of given inner product $\alpha$, which gives the second moment by optimizing over $\alpha$. In Section \ref{sec: S ineq} we compare aforementioned limits for non-negative $\alpha$ in Theorem \ref{thm: S ineq}. Finally in Section \ref{sec: thm exp sec mom proof} we combine our insights to prove Theorem \ref{thm: exp sec mom}.

\section{Formulas for annealed one and two point complexities}\label{sec: moments}
In this section we control the expected exponential rate of the number of local maxima and number of pairs of local maxima (see Lemmata \ref{lem: first mom} and \ref{lem: sec mom} below) by adapting results of \cite[Theorem 5]{arous2020geometry} from critical points to local maxima. We first state our results and then continue with the needed proofs. 
\begin{lemma}[First moment]\label{lem: first mom}
For any open $D\subset\R$ with $2\sqrt{\xipp} \not\in \partial D$
$$ \lim_{N\rightarrow \infty}\frac{1}{N}\ln E[ \Num(D)] = \sup_{x\in D} I(x).$$
\end{lemma}
\noindent To state the result on pairs of local maxima properly consider
$$ \Num_2( D, A ) = \# \left\{(\sigma_1,\sigma_2): \sigma_1\sigma_2\in A,\; \sigma_i \text{ loc. max. of } H_N \text{ with } \frac{1}{N}\partial_r H_N(\sigma_i) \in D \text{ for } i\in\{1,2\} \right\} ,$$
where ``loc. max.'' refers to local maxima on the sphere.  With $\xipa= \xi'(\alpha)$ and $\xippa=\xi''(\alpha)$ as shorthands we define 
\beq \label{def: two-point-function} S(x,\alpha) =\frac{1}{2}\ln\left(\frac{(1-\alpha^2)\xipp^2}{\xip^2-\xipa^2}\right)+\frac{x^2}{2\xipp} Q(\alpha)+2\; \Omega\left(\frac{x}{\sqrt{2\xipp}}\right),\eeq
where 
\begin{equation}\label{eq: Q def}
Q(\alpha)= 1-\frac{2\xipp(\xip-\alpha\xipa+(1-\alpha^2)\xippa)}{\xip^2-\xipa^2+(\xip-\alpha\xipa)(\xipp+\xippa)+(1-\alpha^2)\xipp\xippa}
\end{equation}
 which allows us to state our result on pairs of local maxima as follows: 
\begin{lemma}[Two point annealed complexity at exponential scale]\label{lem: sec mom}
It holds for any interval $A \subset (-1,1)$, that
$$ \lim_{\varepsilon \searrow 0}  \lim_{N\to \infty} \frac{1}{N} \ln E\left[ \Num_2\left( (x-\varepsilon,x+\varepsilon), A\right) \right] \leq \sup_{\alpha\in A} S(x,\alpha).$$
\end{lemma}
The remainder of the section is devoted to the proofs:
\begin{proof}[\bf Proof of Lemma \ref{lem: first mom}]
We adapt the proof of \cite[Theorem 5]{arous2020geometry} with $q=1$, $B=\R$ and consider local maxima instead of critical points. The change to local maxima only causes an extra indicator of the event $\{\lambda_{\max}\left(\nabla^2_{\sph}H_N(\sigma)\right)<0\}$ to appear, when using Kac-Rice formula \cite[Theorem 12.1.1]{adler2009random}, since a critical point $\sigma$ on this event is a local maximum and vice versa almost surely. Following the proof up to \cite[(4.2)]{arous2020geometry}, conditioning on the radial derivative, but not on the energy, and carrying the indicator along we obtain
\beq\label{fstmomeq1} \bea E[ \Num(D)]=& \exp\left(N\left(\frac{1}{2}+\frac{1}{2}\ln\left(\frac{\xipp}{\xip}\right)\right)+o(N)\right) \times \\
 \times& \int_{D} \exp\left(-N\frac{ x^2}{2(\xip+\xipp)}\right) \E\left[\left|\det\left(G-\sqrt{\frac{N}{N-1}}\frac{x}{\xipp}I_{N-1}\right)\right|\mathbf{1}_{\left\{\lambda_{\max}\left(G\right)<\sqrt{\frac{N}{N-1}}\frac{x}{\sqrt{\xipp}}\right\}}\right]dx  ,\eea\eeq
where $I_{N-1}$ is the identity matrix of dimension $N-1$, $\lambda_{\max}$ is the largest eigenvalue and $G$ is a normalized GOE matrix of dimension $N-1$, i.e it is real symmetric with otherwise independent centered Gaussian entries of variance $\frac{1}{N-1}$ off the diagonal and variance $\frac{2}{N-1}$ on the diagonal. 
This reduces the problem to computing the $\frac{1}{N}\ln$ limit of the expectation and applying Laplace principle. Let $\lambda_1\leq \lambda_2\leq ... \leq \lambda_{N-1}$ be the eigenvalues of $G-\sqrt{\frac{N}{N-1}}\frac{x}{\sqrt{\xipp}}I_{N-1}$. If $x<2\sqrt{\xipp}$ we have by Cauchy-Schwartz and estimating roughly
$$ \E\left[\prod_{i=1}^{N-1}|\lambda_i|\mathbf{1}_{\left\{\lambda_{N-1}<0\right\}}\right]\leq \E[\lambda_{N-1}^{2(N-1)}]^{1/2}\P(\lambda_{N-1}<0)^{1/2}.$$
Now by \cite[Lemma 6.3]{arous2001aging} the first term is at most of order $\exp(\const N)$, while the second term decays as $\exp(-\const N^2)$, since it requires a macroscopic change of the empirical spectral distribution of a GOE matrix which has LDP of rate $N^2$ see e.g. \cite[Theorem 6.1]{arous2001aging}. Hence $x<2\sqrt{\xipp}$ contribute $-\infty$ to the $\frac{1}{N}\ln$ limit of the expectation. On the other hand if $x>2\sqrt{\xipp}$ we have 
$$ \E\left[\prod_{i=1}^{N-1}|\lambda_i|\mathbf{1}_{\left\{\lambda_{N-1}<0\right\}}\right]= \E\left[\prod_{i=1}^{N-1}|\lambda_i|\right]-\E\left[\prod_{i=1}^{N-1}|\lambda_i|\mathbf{1}_{\left\{\lambda_{N-1}\geq 0\right\}}\right],$$
where by \cite[(4.3)]{arous2020geometry} we have
\beq \lim_{\epsilon\searrow 0} \lim_{N\rightarrow\infty}\frac{1}{N}\ln\left(\int_{y-\epsilon}^{y+\epsilon}\E\left[\prod_{i=1}^{N-1}|\lambda_i|\right] dx\right) = \frac{y^2}{4\xipp}-\frac{1}{2}-\frac{y}{4\sqrt{\xipp}}\sqrt{\frac{y^2}{\xipp}-4}+\ln\left(\sqrt{\frac{y^2}{4\xipp}-1}+\frac{y}{2\sqrt{\xipp}}\right)\eeq
and by Cauchy-Schwarz followed by application of \cite[Corollary 22 and 23]{subag2017complexity} gives us 
$$ \E\left[\prod_{i=1}^{N-1}|\lambda_i|\mathbf{1}_{\left\{\lambda_{N-1}\geq 0\right\}}\right]\leq \E\left[\prod_{i=1}^{N-1}|\lambda_i|^2\right]^{1/2} \P\left(\lambda_{N-1}\geq 0\right)^{1/2}\leq C\;\E\left[\prod_{i=1}^{N-1}|\lambda_i|\right] \P\left(\lambda_{N-1}\geq 0\right)^{1/2}$$
for some constant $C>0$ independent of $N$. For $x>2\sqrt{\xipp}$ we have by \cite[Theorem 6.2]{arous2001aging} that $\mathbb{P}(\lambda_{N-1}\geq 0)\to0$ and therefore the indicator does not contribute to the limit for $x>2\sqrt{\xipp}$.
Collecting cases we have shown that 
\beq\label{fstmomdetterm}\bea\lim_{\epsilon\searrow 0} \lim_{N\rightarrow\infty}\frac{1}{N}\ln&\left(\int_{y-\epsilon}^{y+\epsilon}\E\left[\prod_{i=1}^{N-1}|\lambda_i|\mathbf{1}_{\left\{\lambda_{N-1}<0\right\}}\right]dx \right) =\\ =&\begin{cases}\frac{y^2}{4\xipp}-\frac{1}{2}-\frac{y}{4\sqrt{\xipp}}\sqrt{\frac{y^2}{\xipp}-4}+\ln\left(\sqrt{\frac{y^2}{4\xipp}-1}+\frac{y}{2\sqrt{\xipp}}\right) & \mbox{ for } y\geq 2\sqrt{\xipp}\\
-\infty & \mbox{ else }\end{cases}\eea\eeq
Applying Laplace principle to \eqref{fstmomeq1} using \eqref{fstmomdetterm} yields the claim.
\end{proof}

\begin{proof}[\bf Proof of Lemma \ref{lem: sec mom}]
Since all local maxima are critical points we have for by \cite[Theorem 6]{arous2020geometry} for $q_1=q_2=1$, $B_{1}=B_{2}=\R$, $D_1=D_2=(x-\varepsilon,x+\varepsilon)$ and $I=A$ that 
$$ \lim_{\varepsilon\searrow0}\limsup_{N\rightarrow\infty}\E[\mathcal{N}_{2}(D_1,I)]\leq \sup_{\alpha\in A, u_1,u_2\in \R} \Psi_{\xi,1,1}(\alpha,u_1,u_2,x,x),$$ 
where $\Psi$ is given in \cite[(3.7)]{arous2020geometry}. 
The optimization with respect to $u_1,u_2$, using Laplace principle on correlated Gaussian tails, yields 
$$ \sup_{u_1,u_2\in \R}-\frac{1}{2}(u_1,u_2,x,x)\Sigma_{U,X}^ {-1}(\alpha,1,1)(u_1,u_2,x,x)^T = -\frac{1}{2}(x,x)\Sigma_{X}^{-1}(\alpha,1,1)(x,x)^T,$$
where by \cite[(A.3)]{arous2020geometry}
$$\Sigma_X(\alpha,1,1)= 
\begin{pmatrix} 
\xipp+\xip- \frac{(\alpha \xippa +\xipa )^2 (1-\alpha^2)\xip}{\xip^2-(\alpha\xipa-\xippa(1-\alpha^2))^2} & \alpha^2 \xippa +\alpha \xipa -\frac{(\alpha \xippa +\xipa )^2 (1-\alpha^2)(\alpha\xipa-\xippa(1-\alpha^2))}{\xip^2-(\alpha\xipa-\xippa(1-\alpha^2))^2} \\ 
\alpha^2 \xippa +\alpha \xipa -\frac{(\alpha \xippa +\xipa )^2 (1-\alpha^2)(\alpha\xipa-\xippa(1-\alpha^2))}{\xip^2-(\alpha\xipa-\xippa(1-\alpha^2))^2}& \xipp+\xip- \frac{(\alpha \xippa +\xipa )^2 (1-\alpha^2)\xip}{\xip^2-(\alpha\xipa-\xippa(1-\alpha^2))^2}
\end{pmatrix}.$$
For $x\geq 2\sqrt{\xipp}$ we obtain the claim, being careful to not confuse $\Omega$ in \eqref{eq: omegadef} with  \cite[(3.3)]{arous2020geometry}, by verifying that 
$$ -\frac{1}{2}(x,x)\Sigma_{X}^{-1}(\alpha,1,1)(x,x)^T +\frac{x^2}{2\xipp} = \frac{x^2}{2\xipp}Q(\alpha),$$
which is a straightforward computation using that $(x,x)\Sigma_{X}^{-1}(\alpha,1,1)(x,x)^T = \frac{2x^2}{\Sigma_{X}(\alpha,1,1)_{1,1}+\Sigma_{X}(\alpha,1,1)_{1,2}}$ and canceling $\xip+\alpha\xipa-\xippa(1-\alpha^2)$.
It remains to show for $x<2\sqrt{\xipp}$ that  
$$ \lim_{\varepsilon\searrow0}\limsup_{N\rightarrow\infty}\E[\mathcal{N}_{2}(D_1,I)] = - \infty.$$ To this end we follow the proof of \cite[Theorem 6]{arous2020geometry} for $q_1=q_2=1$, $B_{1}=B_{2}=\R$, $D_1=D_2=(x-\varepsilon,x+\varepsilon)$ and insert an additional indicator (with the same argument as in the proof of Lemma \ref{lem: first mom}) to adjust for the change from critical points to local maxima. We arrive at \cite[(4.12)]{arous2020geometry}, which reads for our case: 
$$\E[\mathcal{N}_2(D,A)]=C_N\int_{A} \mathcal{D}(\alpha)^{N-3} \mathcal{F}(\alpha) \E\left[\prod_{i=1}^2\left|\det\left(M^{(i)}_{N-1}(\alpha)\right)\right|\mathbf{1}_{E_i}\right]  ,$$
where 
$$C_N = \omega_N \omega_{N-1} \left(\frac{1}{2\pi}(N-1) \frac{\xipp}{\xip}\right)^{N-1} ,\quad \omega_N = 2\frac{\pi^{N/2}}{\Gamma(N/2)},$$
$$ \mathcal{D}(\alpha) = \sqrt{1-\alpha^2} \sqrt{1-\frac{\xipa^2}{\xip^2}}, $$
$$ \mathcal{F}(\alpha) =\sqrt{1-\frac{\xipa^2}{\xip^2}}\sqrt{1-\left(\frac{\alpha\xipa -\xippa(1-\alpha^2)}{\xip}\right)^2}, $$
$$E_i = \{X_i(\alpha) \in \sqrt{N} D, \lambda_{\max}(M^{(i)}_{N-1}(\alpha))<0\},$$
and $X_1(\alpha),X_2(\alpha),M^{(1)}_{N-1}(\alpha),M^{(2)}_{N-1}(\alpha)$ have joint distribution given by \cite[Lemma 15]{arous2020geometry}. The terms $C_N, \mathcal{D}(\alpha)^{N-1}, \mathcal{F}(\alpha)$ are all bounded by $\exp( \const N)$ and therefore the claim follows immediately from the Cauchy-Schwartz inequality if we show that 
$$\E\left[\prod_{i=1}^2\det\left(M^{(i)}_{N-1}(\alpha)\right)^2\right] \leq \exp(\const N),$$ 
since by the eigenvalue interlacing theorem $E_1,E_2$ require a large deviation of the empirical spectral measure of a GOE and therefore by \cite[Theorem 6.1]{arous2001aging} the probabilities of $E_1,E_2$ vanish as $\exp(-\const N^2)$. By Cauchy Schwartz using that $M^{(1)}$ and $M^{(2)}$ have the same distribution and roughly estimating we see that 
\beq\label{eq:rough est} \E\left[\prod_{i=1}^2\det\left(M^{(i)}_{N-1}(\alpha)\right)^2\right]\leq \E\left[\left(|\lambda|_{\max}(M^{(1)}(\alpha))\right)^{4(N-1)}\right].\eeq That this is bounded by $\exp( \const N) $ follows immediately from the tail estimate   
\beq \label{eq:tail estimate}\P \left(|\lambda|_{\max}(M^{(1)}(\alpha))>t\right)\leq\exp(-\const t^2 N).\eeq
This is easily obtained from the fact that $M^{(1)}(\alpha)$ is defined as 
$$ \begin{pmatrix}G_{N-2} & \vec{0}\\ \vec{0}^T & 0 \end{pmatrix}+\begin{pmatrix}0_{N-2} & Z^{(1)}(\alpha)\\ Z^{(1)}(\alpha)^T & 0 \end{pmatrix}+\begin{pmatrix}0 & \vec{0}\\ \vec{0}^T &  Q^{(i)}(\alpha)+\sqrt{\frac{N}{N-1 \xipp}}m_i(\alpha,1,1) \end{pmatrix}-\sqrt{\frac{1}{(N-1)\xipp}}X_1(\alpha) I_{N-1},$$ 
and that all terms have operator norm with tails satisfying  \eqref{eq:tail estimate}. Here the top left entry always has dimension $(N-2)\times(N-2)$ and the bottom right $1\times 1$. Also $G_{N-2}$ a normalized GOE matrix. Lastly $I_{N-1}$ is the dimension $N-1$ identity matrix and $0_{N-2}$, $\vec{0}$ are the zero matrix and vector respectively. For more details see \cite[(4.9) and (4.13)]{arous2001aging}. 
The first term is a normalized GOE matrix and the claimed tail estimate is given by e.g. \cite[Lemma 6.3]{arous2001aging}. The second term has operator norm which is up to multiplicative constant the root of a $\chi^2(N-1)$ distributed random variable. Since the rate function of the $\chi^2(1)$ distribution is asymptotically linear at $\infty$ the claimed tail bound follows. The operator norms of the third and fourth terms are the absolute value of a Gaussian with variance of order $N^{-1}$, which have tails as claimed, which  concludes the argument that \eqref{eq:rough est} is at most $\exp\left(\const N\right)$ and therefore yields the claim. 
\end{proof}

\section{Maximization of two point annealed complexity formula}\label{sec: S ineq}

In this section we prove that the two point complexity $S$ on $[0,1]$ is maximized at $\alpha=0$. 
\begin{theorem}[Second moment complexity formula is maximized at $\alpha=0$]\label{thm: S ineq}
For all $x\in \R$ and all mixtures not of form $\xi(x)=c x^2$ we have 
$$ I(x)\geq 0 \Rightarrow  \sup_{\alpha\in[0,1]} S(x,\alpha) = S(x,0) = 2I(x).$$
\end{theorem}
\begin{rem}

a) We will later in the proof of Theorem \ref{thm: exp sec mom} use simple and general geometric considerations to show that $\alpha<0$ can be ignored, see \eqref{eq: simple geometry} and Remark \ref{rem: easy geometry} below for details.

b) The fact that $\alpha=0$ maximizes $S$ means that most pairs of local maxima of given radial derivative are approximately orthogonal. In the terminology of spin glasses one can say that they are ``replica symmetric''.

\end{rem}
The proof we give here is elementary in nature but by no means trivial. It is guided by seeking a sense of algebraic beauty and could not have been reasonably derived without use of numerical checks and computer algebra. While checking step by step is elementary it is hard even for the authors to see a guiding principle.

\begin{proof}[\bf Proof of Theorem \ref{thm: S ineq}] Let $\alpha\in(0,1)$ and $\xi$ be a mixture not of form $\xi(x)=c x^2$ and consider 
\begin{equation}\label{eq: first step}
    S(x,0)-S(x,\alpha)= \frac{1}{2}\ln\left(\frac{\xip^2-\xipa^2}{(1-\alpha^2)\xip^2}\right)+\frac{x^2}{2\xipp}(Q(0)-Q(\alpha)).
\end{equation}
By the definition \eqref{eq: Q def} of $Q$ we have
$$ Q(0)-Q(\alpha)=
2\xipp\left( \frac{\xip-\alpha\xipa+(1-\alpha^2)\xippa}{\xip^2-\xipa^2+(\xip-\alpha\xipa)(\xipp+\xippa)+(1-\alpha^2)\xipp\xippa}-\frac{1}{\xip+\xipp}\right).$$
Since
$$\begin{array}{l}
\left(\xip+\xipp\right)\left(\xip-\alpha\xipa+(1-\alpha^{2})\xipa'\right)-\left( \xip^{2}-\xipa^{2}+(\xip-\alpha\xipa)(\xipp+\xippa)+(1-\alpha^{2})\xipp\xippa\right)\\=-\alpha\xipa\xip-\xip\alpha^{2}\xippa+\xipa^{2}+\alpha\xipa\xippa\\=\left(\xipa-\alpha\xip\right)\left(\xipa+\alpha\xippa\right)
\end{array},$$
we obtain
\beq\label{eq: Q diff} Q(0)-Q(\alpha) =
\frac{2\xipp}{\xip + \xipp} 
\frac{2\xipp(\xipa-\alpha \xip)(\xipa+\alpha \xippa)}{\xip^2-\xipa^2+(\xip-\alpha\xipa)(\xipp+\xippa)+(1-\alpha^2)\xipp\xippa}<0, \eeq
where the negativity is due to $\xipa-\alpha \xip$ being the only negative term. The next lemma will be used to bound the term $\frac{x^2}{2\xip^2}$ in \eqref{eq: first step}.
\begin{lemma}\label{lem: y bound}
For any mixture not of form $\xi(x)=c x^2$ we have 
$$ I\left(x\right)\geq 0 \Rightarrow 2\leq \frac{x^2}{2\xipp} < \frac{\xipp+\xip}{\xipp-\xip}\ln\left(\frac{\xipp}{\xip}\right).$$
\end{lemma}
\begin{proof}
The lower bound on $\frac{x^2}{2\xipp}$ follows immediately from the definition of $I$, which is $-\infty$ if $\frac{x^2}{2\xipp}<2$. Assume $\frac{x^2}{2\xipp} \geq \frac{\xipp+\xip}{\xipp-\xip}\ln\left(\frac{\xipp}{\xip}\right)$. We will now prove that this implies $I(x)<0$. Note that $h$ given by  
$$ h: [1,\infty)\rightarrow \R \quad x\mapsto \begin{cases} \frac{x+1}{x-1}\ln(x) & \mbox{ for } x>1, \\ 2 & \mbox{ for } x=1,\end{cases}$$
is continuous and strictly increasing. Using that the mixture is not of form $c x^2 $ we have $\xipp>\xip$ and therefore 
$$ \frac{x^2}{2\xipp} \geq \frac{\xipp+\xip}{\xipp-\xip}\ln\left(\frac{\xipp}{\xip}\right) = h\left(\frac{\xipp}{\xip}\right)>h(1) = 2.$$
Using that $\Omega'(y) = -\sqrt{y^2-2}$ and $\Omega(\sqrt{2})=0$ we have $\Omega(y)<0$ for $y>\sqrt{2}$, which gives 
$$ I(x) = \frac{1}{2}\ln\left(\frac{\xipp}{\xip}\right)-\frac{x^2}{4\xipp}\frac{\xipp-\xip}{\xipp+\xip}+\Omega\left(\frac{x}{\sqrt{2\xipp}}\right)
< \frac{1}{2}\ln\left(\frac{\xipp}{\xip}\right)-\frac{\xipp+\xip}{\xipp-\xip}\ln\left(\frac{\xipp}{\xip}\right)\frac{1}{2}\frac{\xipp-\xip}{\xipp+\xip}=0.$$
\end{proof}

Applying Lemma \ref{lem: y bound} to \eqref{eq: first step} while recalling \eqref{eq: Q diff} we have 
\beq\label{eq:S ineq step 1} S(x,0)-S(x,\alpha)> \frac{1}{2}\ln\left(\frac{\xip^2-\xipa^2}{(1-\alpha^2)\xip^2}\right)+\frac{\xipp+\xip}{\xipp-\xip}\ln\left(\frac{\xipp}{\xip}\right)(Q(0)-Q(\alpha)).\eeq
The next lemma bounds the quantity in the $\log$ in this expression.

\begin{lemma}\label{lem: elementary est}
For $\alpha\in[0,1)$ we have 
$$ 1\leq \frac{\xip^2-\xipa^2}{(1-\alpha^2)\xip^2}\leq \frac{\xipp}{\xip}.$$
\end{lemma}
\begin{proof}
For $\alpha=0$ the left inequality is sharp and for $\alpha\rightarrow 1$ the right. Hence the proof is completed by proving that the function 
$$g(\alpha):=\frac{\xip^2-\xipa^2}{(1-\alpha^2)\xip^2}$$
is increasing. By computing $g'$ 
$$ g'(\alpha) = 2 \frac{\alpha(\xip^2-\xipa^2)-(1-\alpha^2)\xipa\xippa}{(1-\alpha^2)^2 \xip^2}$$
we observe that $g$ is increasing if and only if 
\begin{equation}\label{eq: suff to show}
    \frac{\xip^2-\xipa^2}{1-\alpha^2} \geq \frac{\xipa}{\alpha}\xippa.
\end{equation} 
By using that $\xippa$ is non negative and increasing for $\alpha\in[0,1]$ we obtain $\xip-\xipa = \int_{\alpha}^{1}\xippa d\alpha \geq (1-\alpha)\xippa$ which implies 
$$ \frac{\xip^2-\xipa^2}{1-\alpha^2} \geq \xippa\frac{\xip+\xipa}{1+\alpha}= \xippa \sum_{p\geq 2} a_p p \frac{1+\alpha^{p-1}}{1+\alpha}.$$ 
Using that $1\geq \alpha^{p-2}$ we further estimate 
$$ \geq \xippa \sum_{p\geq 2} a_p p \frac{\alpha^{p-2}+\alpha^{p-1}}{1+\alpha} = \xippa \sum_{p\geq 2} a_p p \alpha^{p-2} = \xippa \frac{\xipa}{\alpha},$$
which yields \eqref{eq: suff to show}.
\end{proof}
We can use  Lemma \ref{lem: elementary est} together with the fact that (the continuous continuation on $[1,\infty)$ of) $\frac{\sqrt{z}}{z-1}\ln(z)$ is decreasing together with to obtain 
$$ \frac{\sqrt{\frac{\xip^2-\xipa^2}{(1-\alpha^2)\xip^2}}}{\frac{\xip^2-\xipa^2}{(1-\alpha^2)\xip^2} -1}\ln\left(\frac{\xip^2-\xipa^2}{(1-\alpha^2)\xip^2}\right)\geq \frac{\sqrt{\frac{\xipp}{\xip}}}{\frac{\xipp}{\xip} -1}\ln\left(\frac{\xipp}{\xip}\right). $$
Applying this to \eqref{eq:S ineq step 1} we obtain
$$S(x,0)-S(x,\alpha)> \ln\left(\frac{\xipp}{\xip}\right)\left( \frac{1}{2}\frac{\sqrt{\frac{\xipp}{\xip}}\left(\frac{\xip^2-\xipa^2}{(1-\alpha^2)\xip^2} -1\right)}{\sqrt{\frac{\xip^2-\xipa^2}{(1-\alpha^2)\xip^2}}\left(\frac{\xipp}{\xip} -1\right)}
+\frac{\xipp+\xip}{\xipp-\xip}(Q(0)-Q(\alpha))\right).$$
Since $\ln\left(\frac{\xipp}{\xip}\right)>0$ it suffices to show that 
$$ 4\sqrt{\frac{\xip^2-\xipa^2}{(1-\alpha^2)\xip^2}}\leq \frac{\sqrt{\frac{\xipp}{\xip}}\left(\frac{\xip^2-\xipa^2}{(1-\alpha^2)\xip^2} -1\right)}{\left(\frac{\xipp}{\xip} -1\right)\frac{\xipp+\xip}{\xipp-\xip}\frac{1}{2}(Q(\alpha)-Q(0))}.
$$
Multiplying both sides by $(1-\alpha^2)\sqrt{\frac{\xipp}{\xip}}$ and canceling terms this reads
\beq\label{eq: S ineq step 2} 4\sqrt{1-\left(\frac{\xipa}{\xip}\right)^2}\sqrt{1-\alpha^2}\sqrt{\frac{\xipp}{\xip}}
\leq
\frac{1}{\xip^2}
\frac
{\alpha^2\xip^2-\xipa^2}
{\frac{\xipp+\xip}{\xipp}\frac{1}{2}(Q(\alpha)-Q(0))}.
\eeq
Plugging in the equality in \eqref{eq: Q diff} the right hand side reads 
$$ \frac{\alpha^2-\left(\frac{\xipa}{\xip}\right)^2}{(\alpha \xip-\xipa)(\xipa+\alpha \xippa)}\left(\xip^2-\xipa^2+(\xip-\alpha\xipa)(\xipp+\xippa)+(1-\alpha^2)\xipp\xippa\right).$$
By multiplying the fraction and dividing the bracket by $\xip^2$ and canceling we obtain equality to 
$$ = \frac{\alpha\xip+\xipa}{\xipa+\alpha \xippa}\left(1-\frac{\xipa^2}{\xip^2}+(1-\alpha^2)\frac{\xipp\xippa}{\xip^2}+(1-\alpha\frac{\xipa}{\xip})\frac{\xipp+\xippa}{\xip}\right).$$
Lemma \ref{lem: magic algebra} stated immediately below implies \eqref{eq: S ineq step 2}, by adding \eqref{eq: magic algebra 2} and \eqref{eq: magic algebra 1}. Therefore the proof of Theorem \ref{thm: S ineq} is finished once we have proved Lemma \ref{lem: magic algebra}.

\begin{lemma}\label{lem: magic algebra}
For $\alpha \in (0,1)$, we have
\beq\label{eq: magic algebra 2}2\sqrt{1-\left(\frac{\xipa}{\xip}\right)^2}\sqrt{1-\alpha^2} \sqrt{\frac{\xipp}{\xip}}\leq \left(1-\alpha \frac{\xipa}{\xip}\right) \frac{\xipp+\xippa}{\xip}\frac{\alpha \xip+\xipa}{\xipa+\alpha\xippa} \eeq
as well as 
\beq\label{eq: magic algebra 1} 2\sqrt{1-\left(\frac{\xipa}{\xip}\right)^2}\sqrt{1-\alpha^2} \sqrt{\frac{\xipp}{\xip}}\leq \left(1-\frac{\xipa^2}{\xip^2}+(1-\alpha^2)\frac{\xipp \xippa}{\xip^2}\right)\frac{\alpha \xip+\xipa}{\xipa+\alpha\xippa}.\eeq
\end{lemma}
Inequalities \eqref{eq: magic algebra 2} and \eqref{eq: magic algebra 1} certainly look artificial at first glance, but a surprising structure is in fact hidden under the surface, which is brought to the fore by the following Lemma \ref{lem: reformu}.

The proof of Lemma \ref{lem: magic algebra}
requires certain key observations discussed by Lemmata \ref{lem: reformu}, \ref{lem: SZ X est} and \ref{lem: Z SZ Y est}. Hence we prove these Lemmata first and only then return to the proof of Lemma \ref{lem: magic algebra}.

\begin{lemma}\label{lem: reformu}
With $f(x)= x+\frac{1}{x}$ we have

$$ \eqref{eq: magic algebra 2} \Leftrightarrow 2f(\sqrt{W})\leq f(\sqrt{U})f(\sqrt{V}),$$
$$\mbox{ for } U = \frac{\xipp}{\xippa}, V= \frac{\alpha}{\xipa}\frac{\xip^2-\xipa^2}{(1-\alpha^2)\xip}, W= \frac{\xipa}{\alpha \xippa}.$$
as well as 
$$ \eqref{eq: magic algebra 1} \Leftrightarrow 2f(\sqrt{Z})\leq f(\sqrt{X})f(\sqrt{Y}),$$
$$ \mbox{ for } X= \frac{\xip^2-\xipa^2}{\xipp\xippa(1-\alpha^2)},\mbox{ } Y= \frac{\alpha\xip}{\xipa}, Z= \frac{\xipa}{\alpha\xippa}.$$
\end{lemma}
\begin{proof} 
We start deriving the first equivalence by writing out \eqref{eq: magic algebra 2}, replacing all instances of $\xipp$ by $U\xippa$, as well as dividing by $\sqrt{U}$, which yields 
$$2\sqrt{1-\left(\frac{\xipa}{\xip}\right)^2}\sqrt{1-\alpha^2} \sqrt{\frac{\xippa}{\xip}}\leq \left(1-\alpha \frac{\xipa}{\xip}\right) \frac{f(\sqrt{U})\xippa(\alpha \xip+\xipa)}{\xip(\xipa+\alpha \xippa)}.$$
Next we replace all instances of $\xippa$ by $\frac{\xipa}{\alpha W}$; multiply by $W+1$ and use 
$\left(1-\alpha\frac{\xipa}{\xip}\right)\frac{1}{\alpha\xip}\left(\alpha\xip+\xipa\right)=1-\frac{\xipa^{2}}{\xip^{2}}+\left(1-\alpha^{2}\right)\frac{\xipa}{\alpha\xip}$ to obtain
$$ 2 f(\sqrt{W})\sqrt{1-\left(\frac{\xipa}{\xip}\right)^2}\sqrt{1-\alpha^2} \sqrt{\frac{\xipa}{\alpha\xip}}\leq f(\sqrt{U})\left(1-\frac{\xipa^2}{\xip^2}+ (1-\alpha^2)\frac{\xipa}{\alpha\xip}\right).$$
Dividing by the roots on the left hand side and using the definition of $V$ we easily obtain equivalence to 
$$ 2f(\sqrt{W})\leq f(\sqrt{U})f(\sqrt{V})$$ 
as claimed.

For the second equivalence we proceed in the same fashion. Writing down \eqref{eq: magic algebra 1}, replacing all instances of $\xipp$ with $\frac{\xip^2-\xipa^2}{X\xippa(1-\alpha^2)}$ and multiplying by $\sqrt{X}$ we have 
$$2 \left(1-\left(\frac{\xipa}{\xip}\right)^2\right) \sqrt{\frac{\xip}{\xippa}}\leq f(\sqrt{X})\left(1-\frac{\xipa^2}{\xip^2}\right)\frac{\alpha \xip+\xipa}{\xipa+\alpha\xippa}.$$
After removing the factor $1-\frac{\xipa^2}{\xip^2}$ from both sides we replace all remaining instances of $\xippa$ with $\frac{\xipa}{\alpha Z}$ and multiply by $1+\frac{1}{Z}$ to obtain
$$2 f(Z) \sqrt{\frac{\alpha \xip}{\xipa}}\leq f(\sqrt{X})\frac{\alpha \xip+\xipa}{\xipa}.$$
Reorganizing the remaining terms and replacing with $Y$ we immediately obtain the claim.
\end{proof}

The first part of Lemma \ref{lem: magic algebra} will follow easily from the next lemma.
\begin{lemma}\label{lem: VW ineq}
For all $\alpha\in(0,1)$ it holds that $$V\geq W^{-1}\ge1,$$
for $V,W$ as in Lemma \ref{lem: reformu}.
\end{lemma}
\begin{proof}
Since we only consider $p\geq 2$ it is easy to check that $W^{-1}\geq 1$. The claim $V\geq W^{-1}$ is equivalent to
\beq\label{eq: almost done 1} \frac{\xipa^2}{\xip^2}+(1-\alpha^2)\frac{\xippa}{\xip}\leq 1.\eeq
Consider the distribution given by $\mathbb{P}(P=p) \propto a_p p$ and write \eqref{eq: almost done 1} as follows
$$ \E[\alpha^{P-1}]^2+(1-\alpha^2)\E[(P-1)\alpha^{P-2}]\leq 1.$$
This inequality follows by estimating $\E_Q[\alpha^{P-1}]^2\leq \E_Q[\alpha^{2P-2}]$
and realizing that 
\beq \label{eq: almost done 2} \alpha^{2P-2}+(1-\alpha^2)(P-1)\alpha^{P-2}\leq 1 \eeq 
for all $\alpha\in [0,1]$ and $P\geq 2$.

To verify \eqref{eq: almost done 2} first note for $\alpha \in[0,1]$
\begin{equation}\label{eq: deriv}
    \frac{\partial}{\partial \alpha}\left( \alpha^{2P-2}+(1-\alpha^2)(P-1)\alpha^{P-2}\right)
    =
    (P-1) \alpha^{P-3} \left( 2\alpha^{P}  - 2\alpha^2 + (P-2)(1-\alpha^2)\right).
\end{equation}    
Letting
$$ h_P(\alpha):= 2\alpha^{P} - 2\alpha^2 + 2\alpha^2 + (P-2)(1-\alpha^2)=2\alpha^{P} + P-2 - P \alpha^2,$$
we have
$$h_P'(\alpha) =2 P (\alpha^{P-1} - \alpha)\leq 0, $$
and therefore $h_P(\alpha)\geq h_P(1)=0$. From this we see that \eqref{eq: deriv} is non-negative, so \eqref{eq: almost done 1} follows since it trivially holds for $\alpha=1$.
\end{proof}

We are now ready to prove the first part of Lemma \ref{lem: magic algebra}.
\begin{proof}[\bf Proof of \eqref{eq: magic algebra 2}]
By Lemma \ref{lem: reformu} and using its notation it suffices to show that
$$ 2f(\sqrt{W})\leq f(\sqrt{U})f(\sqrt{V}).$$ 
Clearly $f$ is always at least $2$ and therefore this in turns follows from
$$ f(\sqrt{W})\leq f(\sqrt{V}),$$
which is a consequence of Lemma \ref{lem: VW ineq}. This completes the proof of \eqref{eq: magic algebra 2}.
\end{proof}

Before proving the second part of Lemma \ref{lem: magic algebra} we prove two inequalities involving
$$ S= \frac{Y+1}{2}\left(1-\frac{\alpha}{2}\frac{Y-1}{Y}\right).$$
The first is:
\begin{lemma}\label{lem: Z SZ Y est} 
We have for $\alpha\in(0,1)$
\beq\label{eq: Z SZ Y est}2f(\sqrt{Z})\leq f(\sqrt{S Z})f(\sqrt{Y}).\eeq
\end{lemma}
\begin{proof}
Dividing \eqref{eq: Z SZ Y est} by $\sqrt{Z}$ we obtain the equivalent representation 
$$ 2+2 Z^{-1} \leq \sqrt{S}f(\sqrt{Y})+ Z^{-1}\frac{f(\sqrt{Y})}{\sqrt{S}}.$$
Hence it suffices to show 
$$ Z^{-1}\left(2-\frac{f(\sqrt{Y})}{\sqrt{S}}\right)\leq \sqrt{S}f(\sqrt{Y})-2.$$
By observing that
\beq\label{eq: almost done 4} \frac{f(\sqrt{Y})^2}{4 S}= \frac{Y+2+Y^{-1}}{2(Y+1)\left(1-\frac{\alpha}{2}\frac{Y-1}{Y}\right)} = 1- \frac{(1-\alpha)(Y-Y^{-1})}{2(Y+1)\left(1-\frac{\alpha}{2}\frac{Y-1}{Y}\right)}< 1,\eeq
since $\alpha\in(0,1)$ and $Y > 1$, we obtain $2-\frac{f(\sqrt{Y})}{\sqrt{S}}\geq 0$. Therefore \eqref{eq: Z SZ Y est} is also equivalent to 
\beq\label{eq: almost done 6}Z^{-1}\leq \frac{\sqrt{S}f(\sqrt{Y})-2}{2-\frac{f(\sqrt{Y})}{\sqrt{S}}}.\eeq
Reformulating the right hand side of this inequality yields 
$$ \frac{\sqrt{S}f(\sqrt{Y})-2}{2-\frac{f(\sqrt{Y})}{\sqrt{S}}} =  \frac{(\sqrt{S}f(\sqrt{Y})-2)(2+\frac{f(\sqrt{Y})}{\sqrt{S}})}{4-\frac{f(\sqrt{Y})^2}{S}} = \frac{S(f(\sqrt{Y})^2-4)+2f(\sqrt{Y})(S-1)\sqrt{S} }{4 S-f(\sqrt{Y})^2}.$$
\beq\label{eq: almost done 5} = \frac{
\left\{(S-1) f(\sqrt{Y})^2+f(\sqrt{Y})^2 -4 S\right\}+(S-1)f(\sqrt{Y})^2\sqrt{\frac{4 S}{f(\sqrt{Y})^2} }}{4 S-f(\sqrt{Y})^2}. \eeq
We then use the easily checked representations 
$$ S= \frac{(Y+1)^2}{4 Y}\left(1+(1-\alpha)\frac{Y-1}{Y+1} \right) \mbox{ and } f(\sqrt{Y})^2=\frac{(Y+1)^2}{Y}$$
to compute \eqref{eq: almost done 5} piece by piece as follows
$$ S-1 = \frac{Y-1}{4 Y} \left( (2-\alpha)Y -\alpha\right),$$
$$ (S-1)f(\sqrt{Y})^2 = \left( Y(2-\alpha)-\alpha\right) \frac{(Y-1)(Y+1)^2}{4 Y^2},$$
$$ \frac{4 S}{f(\sqrt{Y})^2} =1+(1-\alpha)\frac{Y-1}{Y+1},$$
$$ 4S-f(\sqrt{Y})^2 = (1-\alpha)\frac{(Y-1)(Y+1)}{Y}.$$ 
Using the computations so far we also obtain for the curly bracket in \eqref{eq: almost done 5}
$$ (S-1)f(\sqrt{Y})^2-4S+f(\sqrt{Y})^2= \left(\left( Y(2-\alpha)-\alpha)\right) (Y+1)-4(1-\alpha)Y\right)\frac{(Y-1)(Y+1)}{4Y^2},$$
where the first factor on the LHS equals
$$
\begin{array}{l}
\left(Y\left(1+1-\alpha\right)-1+1-\alpha\right)\left(Y+1\right)-4\left(1-\alpha\right)Y\\
=\left(Y-1\right)\left(Y+1\right)+\left(1-\alpha\right)\left\{ \left(Y+1\right)\left\{ Y+1\right\} -4Y\right\} \\
=Y^{2}-1+\left(1-\alpha\right)\left(Y-1\right)^{2}=\left(Y^{2}-1\right)\left(1+\left(1-\alpha\right)\frac{Y-1}{Y+1}\right).
\end{array}$$

Collecting the pieces of \eqref{eq: almost done 5} we computed and multiplying enumerator and denominator by $\frac{4 Y^2}{(Y+1)(Y-1)}$ we obtain 
$$ \frac{\sqrt{S}f(\sqrt{Y})-2}{2-\frac{f(\sqrt{Y})}{\sqrt{S}}}=\frac{ (Y^2-1)\left(1+  (1-\alpha)\frac{Y-1}{Y+1}\right)  + \left( Y(2-\alpha)-\alpha\right) (Y+1)\sqrt{1+(1-\alpha)\frac{Y-1}{Y+1}} }{4(1-\alpha)Y}. $$
By using the trivial estimate $1+(1-\alpha)\frac{Y-1}{Y+1}\geq 1$ (recall $Y>1$) twice we clearly have
$$  \frac{\sqrt{S}f(\sqrt{Y})-2}{2-\frac{f(\sqrt{Y})}{\sqrt{S}}} \geq \frac{ Y^2-1  + \left( Y(2-\alpha)-\alpha\right) (Y+1) }{4(1-\alpha)Y}= \frac{1}{4}\frac{Y+1}{Y}\left(2\frac{Y-1}{1-\alpha}+Y+1\right).$$
Estimating further using $Y+1\geq 2$ then yields 
\beq\label{eq: almost done 7}  \frac{\sqrt{S}f(\sqrt{Y})-2}{2-\frac{f(\sqrt{Y})}{\sqrt{S}}} \geq  \frac{1}{2}\frac{Y+1}{Y}\left(\frac{Y-1}{1-\alpha}+1\right).\eeq

It remains to show that 
\beq\label{eq: almost done 8} Z^{-1} \leq \frac{1}{2}\frac{Y+1}{Y}\left(\frac{Y-1}{1-\alpha}+1\right),\eeq
since then \eqref{eq: almost done 6} and thereby the claim immediately follow from \eqref{eq: almost done 7}. By definition of $Y$ and $Z$, after dividing both sides by $Y$,   \eqref{eq: almost done 8} reads
\beq\label{eq: almost done 9}\frac{\alpha\xippa}{\xip} \leq \frac{1}{2}\frac{\xipa}{\alpha\xip}\left(1+\frac{\xipa}{\alpha\xip}\right)\left(\frac{\frac{\alpha\xip}{\xipa}-1}{1-\alpha}+1\right)= \frac{1}{2}\left(1+\frac{\xipa}{\alpha\xip}\right)\frac{1-\frac{\xipa}{\xip}}{1-\alpha}.\eeq
With $P\geq 2$ a random variable with $\P(P=p)\propto a_p p$ we can write \eqref{eq: almost done 9} as 
\begin{equation}\label{eq: expectation inequality}
    \E[(P-1)\alpha^{P-2}]\leq \frac{1}{2}\left(1+\E[\alpha^{P-2}]\right)\frac{1-\E[\alpha^{P-1}]}{1-\alpha}.
\end{equation}
Multiplying by $2(1-\alpha)$ and bringing all terms except the $1$ to the left hand side we have the equivalent representation 
$$\E[2(P-1)(1-\alpha)\alpha^{P-2}]-(1-\alpha)\E[\alpha^{P-2}]+\alpha\E[\alpha^{P-2}]^2\leq 1.$$
But by Cauchy Schwartz inequality and linearity of expectation
$$ \E[2(P-1)(1-\alpha)\alpha^{P-2}]-(1-\alpha)\E[\alpha^{P-2}]+\alpha\E[\alpha^{P-2}]^2\leq \E[\alpha^{P-2}\left((2P-3)(1-\alpha)+\alpha^{P-1}\right)].$$
Then \eqref{eq: expectation inequality} and \eqref{eq: almost done 9} follows once we have show that
\beq\label{eq: almost done 10} \alpha^{P-2}\left((2P-3)(1-\alpha)+\alpha^{P-1}\right) \leq 1 \text{ for all }\alpha \in [0,1], P\ge 2.\eeq
To this end note that
$$\frac{\partial}{\partial \alpha} \left( \alpha^{P-2}\left(\alpha^{P-1}+(2P-3)(1-\alpha)\right)\right) = (2P-3)\alpha^{P-3}h_P(\alpha)$$
for 
$$ h_P(\alpha):= P-2-(P-1)\alpha+\alpha^{P-1},$$
and $h_P(\alpha)\ge h_P(1)=0$ by checking that $h_P'(\alpha)\le0$. Since \eqref{eq: almost done 10} trivially holds for $\alpha=1$ this proves \eqref{eq: almost done 10}, and finishes the proof. 
\end{proof}

The second inequality we need for the second second part of Lemma \ref{lem: magic algebra} is:

\begin{lemma}\label{lem: SZ X est}
We have for $\alpha\in(0,1)$
$$\xippa \geq \xip \Rightarrow f(\sqrt{S Z}) \le f(\sqrt{X}).$$
\end{lemma}
\begin{proof}
The claim follows immediately from 
$$ X\leq S Z \leq 1. $$
The second inequality is quickly checked by writing it out partially 
$$ S Z = \frac{\alpha\xip+\xipa}{2\alpha\xippa}\left(1-\frac{\alpha}{2}\frac{Y-1}{Y}\right)$$
and observing that $\xip\leq \xippa$ by assumption, $\xipa\leq \alpha \xippa$ trivially, remembering that $Y \geq 1 $ and checking that both factors are at most $1$.

It remains to check $X\leq S Z$. By definition 
$$ X = \frac{\xip^2-\xipa^2}{\xipp\xippa(1-\alpha^2)}=\frac{\int_{\alpha}^{1}  \xi''(t)dt}{1-\alpha} \frac{\xip+\xipa}{\xipp\xippa(1+\alpha)}.$$
Since $\alpha\in[0,1]$ and $\xi''$ is convex on $[0,1]$ we have $\frac{\int_{\alpha}^{1}  \xi''(t)dt}{1-\alpha}\le\frac{\xipp+\xippa}{2}$ and so
$$ X\leq \frac{\xipp+\xippa}{2}\frac{\xip+\xipa}{\xipp\xippa(1+\alpha)}= \frac{\xipa}{2\alpha\xippa}\left(1+\frac{\xippa}{\xipp}\right)\left(1+\frac{\xip}{\xipa}\right) \frac{\alpha}{1+\alpha}.$$
Consider random variables $P,Q\geq 2$ with $\P(P=p)\propto a_p p$ and $\P(Q=p)\propto a_p p(p-1)$. Clearly $\P(P\geq k)\leq \P(Q\geq k)$ and so since $\alpha \in (0,1)$ and $P\ge2$
\beq 
\frac{\xippa}{\xipp}= \E[\alpha^{Q-2}] \leq \E[\alpha^{P-2}] = \frac{\xipa}{\alpha \xip}.\eeq
Thus
$$ X\leq  \frac{\xipa}{2\alpha\xippa}\left(1+\frac{\xipa}{\alpha\xip}\right)\left(1+\frac{\xip}{\xipa}\right) \frac{\alpha}{1-\alpha}= \frac{Z}{2} (1+Y^{-1}) \frac{\alpha+ Y}{1+\alpha} = Z \frac{Y+1}{2}\left(1-\frac{\alpha}{1+\alpha}\frac{Y-1}{Y}\right) .$$
Remembering that $Y\geq 1$ and estimating $1+\alpha\leq 2$ yields $X\leq Z S$, i.e. the claim.
\end{proof}

We are now ready to prove the second part of Lemma \ref{lem: magic algebra}.
\begin{proof}[\bf Proof of \eqref{eq: magic algebra 1}]
By Lemma \ref{lem: reformu} we need to prove 
$$ 2f(\sqrt{Z})\leq f(\sqrt{X})f(\sqrt{Y}). $$
To this end first note that if $\xip\geq \xippa$ the claim follows immediately, since then $Y\geq Z^{-1}\geq 1$ and therefore $f(\sqrt{Y})\geq f(\sqrt{Z})$, using that $f(x)\ge2$. Hence we can assume $\xip<\xippa$ and by Lemma \ref{lem: Z SZ Y est} and Lemma \ref{lem: SZ X est} the claim follows.
\end{proof}
This completes the proof of Lemma \ref{lem: magic algebra}, and as explained above the statement of that lemma, it also  completes the proof of Theorem \ref{thm: S ineq}.
\end{proof}

\section{Proof of the main results}\label{sec: thm exp sec mom proof}
This section is devoted to the proof of Theorem \ref{thm: exp sec mom}.

\begin{proof}[\bf Proof of Theorem \ref{thm: exp sec mom}]
Equation \eqref{eq: exp fst mom} follows immediately from Lemma \ref{lem: first mom} and Markov inequality, recalling \eqref{def: rinfty}. 
When $\xi(x)=c x^2$ one can check that $S(x,\alpha)=  2\Omega(\frac{x}{2\sqrt{c}})$ for all $\alpha$ so the claim trivially follows. Hence we assume the contrary for the remainder of the proof.

To verify \eqref{eq: exp sec mom} note that for any measurable set $D$ by Cauchy Schwartz inequality we have 
$$ \E[\mathcal{N}(D)]^2 \le \E[\mathcal{N}(D)^2],$$
hence it is sufficient to show that the reverse inequality holds on exponential scale. 

To this end note that for any $\delta>0$
\beq\label{eq: simple geometry} 0 \leq \left|\sum_{\sigma\in \mathcal{N}(D)}\sigma\right|_2^2=
\sum_{(\sigma,\tau)\in \mathcal{N}_2(D,[-1,1])} \sigma \tau \leq \mathcal{N}_2(D,(-\delta,1])-\delta \mathcal{N}_2(D,[-1,-\delta])\eeq
and therefore 
\beq \label{eq: neg alph est}\mathcal{N}_2(D,[-1,1])\leq (\frac{1}{\delta}+1) \mathcal{N}_2(D,[-\delta,1]).\eeq
Using $\mathcal{N}_2(D,[-1,1]) = \mathcal{N}(D)^2$ we obtain
\[
\begin{array}{l}
\lim_{\varepsilon\searrow0}\lim_{N\rightarrow\infty}\frac{1}{N}\ln\E[\mathcal{N}((x-\varepsilon,x+\varepsilon))^{2}]\\
\le\lim_{\varepsilon\searrow0}\lim_{N\rightarrow\infty}\frac{1}{N}\ln\E\left[\mathcal{N}\left((x-\varepsilon,x+\varepsilon)\right)+\mathcal{N}_{2}\left((x-\varepsilon,x+\varepsilon),\left[-\delta,1\right)\right)\right].
\end{array}
\]
Applying Lemmata \ref{lem: first mom} and \ref{lem: sec mom} yields that this is at most
$$ \max\left\{\sup_{\alpha\in \left[-\delta,1\right)} S(x,\alpha),I(x)\right\}.  $$
Since $S$ is continuous in $\alpha$ we may optimize over $\delta>0$. Furthermore by assumption $x\in [r_\infty,r_0]$ and therefore $S(x,0)=2 I(x)\geq I(x)\geq 0$. Hence we obtain overall 
$$\lim_{\varepsilon\searrow 0}\lim_{N\rightarrow\infty}\frac{1}{N}\ln\left(\E[\mathcal{N}((x-\varepsilon, x+\varepsilon))^2]\right)\leq \sup_{\alpha\in[0,1]} S(x,\alpha).$$
By Theorem \ref{thm: S ineq} this supremum is attained in $\alpha=0$ and the claim follows since
$$ \sup_{\alpha\in[0,1]} S(x,\alpha)= S(x,0)=2 I(x) =\lim_{\varepsilon\searrow 0}\lim_{N\rightarrow\infty}\frac{2}{N}\ln\left(\E[\mathcal{N}((x-\varepsilon, x+\varepsilon))]\right). $$
\end{proof}
\begin{rem}\label{rem: easy geometry}
The estimate \eqref{eq: neg alph est} can be improved to
$$ \E[\mathcal{N}_2(D,[-1,1])]\leq (1 +o(1)) \E[\mathcal{N}_2(D,[-\delta,1])]$$ 
by replacing \eqref{eq: simple geometry} with 
$$ 0 \leq \left|\sum_{\sigma\in \mathcal{N}(D)}\sigma\right|_2^2=
\sum_{(\sigma,\tau)\in \mathcal{N}_2(D,[-1,1])} \sigma \tau \leq \delta^2\mathcal{N}_2(D,(-\delta,\delta^2])-\delta \mathcal{N}_2(D,[-1,-\delta])+ \mathcal{N}_2(D,(\delta^2,1])$$
and applying Lemmata \ref{lem: first mom}, \ref{lem: sec mom} and Theorem \ref{thm: S ineq}. While this is not necessary on an exponential scale it will be necessary if one wants to show matching of moments up to multiplicative error $1+o(1)$ to prove concentration on the mean as in \cite{subag2017complexity}.

\end{rem}

\printbibliography
\end{document}